\documentclass[12pt]{amsart}

\usepackage{latexsym, amssymb, amscd, amsfonts}
\usepackage{bm}
\usepackage[colorlinks=true,urlcolor=black,citecolor=black,linkcolor=black,%
pdftitle={Positive systems of Kostant roots},%
pdfauthor={Ivan Dimitrov and Mike Roth},%
pdfsubject={Representation theory},%
pdfkeywords={Parabolic subalgebras, Kostant root systems, Positive roots}]{hyperref}
\usepackage{ifthen}

\DeclareMathAlphabet{\mathpzc}{OT1}{pzc}{m}{it}

\newboolean{bourbaki}
\setboolean{bourbaki}{true}

\newcommand{\np}{\medskip\noindent}

\newcounter{eqcounter}[section]
\renewcommand{\theeqcounter}{\arabic{section}.\arabic{eqcounter}}
\renewenvironment{equation}{\medskip\noindent\refstepcounter{eqcounter}\makebox[0pt][l]{({\bf\theeqcounter})}\begin{minipage}[b]{\textwidth}$$}{$$\end{minipage}\medskip\noindent}

\newcommand{\Old}[1]{\oldstylenums{#1}}

\newcommand{\NumberList}{\renewcommand{\labelenumi}{{\bf \Old{\arabic{enumi}}.}}}
\newcommand{\AlphaList}{\renewcommand{\labelenumi}{{({\em\alph{enumi}})}}}
\newcommand{\RomanList}{\renewcommand{\labelenumi}{{({\roman{enumi}})}}}

\renewcommand{\geq}{\geqslant}
\renewcommand{\leq}{\leqslant}

\newcommand{\st}{\,\,|\,\,}                            % "such that" separator in set descriptions
\newcommand{\Msh}{{\mathcal M}}                        %  Submodule of \Nsh
\renewcommand{\Rsh}{{\mathcal R}}                      %  the set of a-roots
\newcommand{\Ssh}{{\mathcal S}}  				%  Subset of the a-roots 
\newcommand{\Zsh}{{\mathcal Z}}   
\newcommand{\Tsh}{{\mathcal T}}                        

\newcommand{\A}{\mathrm{A}}                            %  A_n root system
\newcommand{\B}{\mathrm{B}}                            %  Borel subgroup, root system of B_n type
\newcommand{\C}{\mathrm{C}}                            %  Root system of C_n type
\newcommand{\D}{\mathrm{D}}                            %  Root system of D_n type
\newcommand{\E}{\mathrm{E}}                            %  Root system of E type
\newcommand{\F}{\mathrm{F}}                            %  Root system of type F
\newcommand{\G}{\mathrm{G}}                            %  Group G
\newcommand{\I}{\mathrm{I}}

\newcommand{\X}{\mathrm{X}}                            %  G/B
\newcommand{\V}{\mathrm{V}}                            %  Variety V, or open set V
\newcommand{\U}{\mathrm{U}}                            %  Open set U
\newcommand{\WW}{\mathrm{W}}                           %  (briefly) a variety
\renewcommand{\SS}{\mathrm{S}}                         %  Subset
\newcommand{\Sym}{\operatorname{Sym}}                  %  Symmetric product

\newcommand{\Ad}{\operatorname{Ad}}
\newcommand{\tr}{\operatorname{tr}}
\newcommand{\supp}{\operatorname{supp}}
\newcommand{\cP}{\mathcal{P}}
\newcommand{\cQ}{\mathcal{Q}}
\def\vep{\varepsilon}

\def\gc{\mathfrak{c}}

\def\gg{\mathfrak{g}}

\def\gh{\mathfrak{h}}

\def\gm{\mathfrak{m}}

\def\gp{\mathfrak{p}}

\def\gs{\mathfrak{s}}
\def\gt{\mathfrak{t}}

\def\ggl{\mathfrak{gl}}
\def\gsl{\mathfrak{sl}}

\newcommand{\opp}{\operatornamewithlimits{\oplus}}

\ifthenelse{\boolean{bourbaki}}
{
\newcommand{\CC}{\mathbf{C}} % complex numbers
\newcommand{\QQ}{\mathbf{Q}} % rational numbers
\newcommand{\RR}{\mathbf{R}} % real numbers
}
{
\newcommand{\CC}{\mathbb{C}} % complex numbers
\newcommand{\QQ}{\mathbb{Q}} % rational numbers
\newcommand{\RR}{\mathbb{R}} % real numbers
}

\newcommand{\und}{\underline}

\newcommand{\point}{\refstepcounter{subsection}\noindent{\bf \thesubsection.} }

\begin{document}
\pagestyle{plain} \title{{ \large{Positive systems of Kostant roots}}
}
\author{Ivan Dimitrov}
\address{Department of Mathematics and Statistics, Queen's University, Kingston,
Ontario,  K7L 3N6, Canada} 
\email{dimitrov@mast.queensu.ca} 
\thanks{Research of I.\ Dimitrov and M.\ Roth was partially supported by NSERC grants}
\author{Mike Roth}
\email{mikeroth@mast.queensu.ca}

\subjclass[2010]{Primary 17B22; Secondary 17B20, 17B25}

\begin{abstract} Let $\gg$ be a simple complex Lie algebra and let $\gt \subset \gg$ be a toral subalgebra of $\gg$. 
As a $\gt$-module $\gg$ decomposes as 
\[\gg = \gs \oplus \big(\oplus_{\nu \in \Rsh} \gg^\nu\big)\]
where $\gs \subset \gg$ is the reductive part of a parabolic subalgebra of $\gg$ and $\Rsh$ is the Kostant root system associated to $\gt$.
When $\gt$ is a Cartan subalgebra of $\gg$ the decomposition above is nothing but the root decomposition of $\gg$ with respect to $\gt$; 
in general the properties of $\Rsh$ resemble the properties of usual root systems. In this note we study the following problem: ``Given a subset 
$\Ssh \subset \Rsh$, is there a parabolic subalgebra $\gp$ of $\gg$ containing $\Msh = \oplus_{\nu \in \Ssh} \gg^\nu$  and whose 
reductive part equals $\gs$?''. Our main results is that, for a classical simple Lie algebra $\gg$ and a saturated $\Ssh \subset \Rsh$, the condition
$(\Sym^\cdot(\Msh))^{\gs} = \CC$ is necessary and sufficient for the existence of such a $\gp$.  In contrast, we show that this statement is no longer 
true for the exceptional
Lie algebras $\F_4, \E_6, \E_7$, and $\E_8$. Finally, we discuss the problem in the case when $\Ssh$ is not saturated.

\np
Keywords: Parabolic subalgebras, Kostant root systems, Positive roots.
\end{abstract}

\maketitle

\vspace{0.5cm}

\section{Introduction}

\np
\point
Let $\gg$ be a simple complex Lie algebra and let $\gh \subset \gg$ be a Cartan subalgebra. 
The root decomposition of $\gg$ with respect to $\gh$ is

\begin{equation}\label{eqn:h-decomp}
\gg = \gh  \oplus \big(\oplus_{\alpha \in \Delta} \gg^\alpha\big)
\end{equation}
where, for any $\alpha \in \gh^*$, 
\[\gg^\alpha:= \{x \in \gg \, | \, [t,x] = \alpha(t) x {\text { for every }} t \in \gh\} \quad \quad \text{and} 
\quad \quad \Delta = \{\alpha \in \gh^*\setminus  \{0\} \, | \, \gg^\alpha \neq 0\}.\]
The Borel subalgebras of $\gg$ containing $\gh$ are in a bijection with the {\it positive systems} $\Delta^+ \subset \Delta$, i.e.,  
the subsets $\Delta^+$ satisfying the following properties: (i) $\Delta = \Delta^+ \cup (-\Delta^+)$, (ii) $\Delta^+ \cap (-\Delta^+) = \emptyset$, and
(iii) $\alpha, \beta \in \Delta^+$, $\alpha + \beta \in \Delta$ implies $\alpha + \beta \in \Delta^+$. Positive systems of roots represent a much 
studied and well-understood topic in the theory of semisimple Lie algebras. Here is a particular problem 
that arises in various situations: ``Given 
a subset $\Phi \subset \Delta$, determine if there is a positive system $\Delta^+$ containing $\Phi$''. 
The answer is that such a positive system 
exists if and only if the semigroup generated by $\Phi$ does not contain $0$. 
The aim of this paper is to address the analogous problem in a more general situation.

\np 
%\bpoint{$\boldsymbol\gt$-roots}
\point
Let $\gt \subset \gg$ be a toral subalgebra of $\gg$, that is, a commutative subalgebra of semisimple elements. 
As a $\gt$-module $\gg$ decomposes as 

\begin{equation}\label{eqn:decomp}
\gg = \gs \oplus \big(\oplus_{\nu \in \Rsh} \gg^\nu\big)
\end{equation}
where 
\[\gg^\nu:= \{x \in \gg \, | \, [t,x] = \nu(t) x {\text { for every }} t \in \gt\}, \quad \gs = \gg^0, \quad \text{and} 
 \quad \Rsh = \{\nu \in \gt^*\setminus  \{0\} \, | \, \gg^\nu \neq 0\}.\] 
 We refer to $\Rsh$ as the  $\gt$-root system of $\gg$, to the elements of $\Rsh$ as  the $\gt$-roots, and
 to the spaces $\gg^\nu$ as the $\gt$-root spaces. 
Often we will drop the reference to $\gt$ when it is clear from the context.

\np
To explain the relation between the decompositions \eqref{eqn:h-decomp} and \eqref{eqn:decomp},
extend $\gt$ to a Cartan subalgebra $\gh$.  The inclusion $\gt\subset\gh$
then induces a surjection $\gh^{*}\rightarrow\gt^{*}$. The $\gt$-root system $\Rsh$ consists of the nonzero elements of the image of $\Delta$ under this
map, and for any $\nu\in \Rsh$ the $\gt$-root space $\gg^{\nu}$ is the sum of the $\gh$-root spaces $\gg^{\alpha}$
such that $\alpha\mapsto\nu$.  Since $\gt$ may be an arbitrary complex subspace of $\gh$ we see that, in contrast
to the case of an $\gh$-decomposition, $\gt$-root spaces may be more than one-dimensional, and $\gt$-roots may be
complex multiples of one another.  (For $\gh$-root systems, $\alpha,r\alpha\in \Delta$ implies that $r=\pm 1$.)

\np
\point
 The subalgebra $\gs$ is a reductive subalgebra of $\gg$ and, moreover, $\gs$ is a reductive part of a parabolic subalgebra of $\gg$.
 Note that $\gt$ is contained in $\Zsh(\gs)$, the centre of $\gs$. In the case when $\gt = \Zsh(\gs)$ the properties of $\Rsh$ and the decomposition 
\eqref{eqn:decomp} were studied by Kostant, \cite{K}. Kostant proved that, for every $\nu \in \Rsh$, $\gg^\nu$ is an irreducible $\gs$-module and showed 
that $\Rsh$ inherits many of the properties of $\Delta$. 
To recognize Kostant's contribution, we refer to the elements of $\Rsh$ 
as ``Kostant roots'' in the title, however we use the shorter ``$\gt$-roots'' in the text. 

\np 
\point
To describe and motivate the problem we address in this note, we assume in this subsection 
that $\gt = \Zsh(\gs)$. We caution the reader that not all of equivalences in the 
following discussion hold when $\gt\neq \Zsh(\gs)$.

\np
One introduces the notion of a positive system $\Rsh^+ \subset \Rsh$ exactly as above: (i) $\Rsh = \Rsh^+ \cup (-\Rsh^+)$, (ii) $\Rsh^+ \cap (-\Rsh^+) = 
\emptyset$, and (iii) $\mu, \nu \in \Rsh^+$, $\mu + \nu \in \Rsh$ implies $\mu + \nu \in \Rsh^+$. Proposition VI.1.7.20 in \cite{Bo} implies
that positive systems in $\Rsh$
are in a bijection with parabolic subalgebras of $\gg$ whose reductive part is $\gs$. The paper 
\cite{DFG} contains a detailed 
discussion (in slightly different terms) of positive systems $\Rsh^+$. In particular, a result of \cite{DFG} implies that a subset $\Tsh \subset \Rsh$ is a
positive system if and only if there exists a linear function $\varphi: \V \to \RR$, $\V$ being the real vector space spanned by $\Rsh$, such that
$\ker \varphi \cap \Tsh = \emptyset$ and $\nu \in \Tsh$ if and only if $\varphi(\nu) > 0$. 
Note that every positive system $\Rsh^+$ is {\it saturated}, i.e., 
$\nu \in \Rsh^+, r \in \QQ_+$ and $r \nu \in \Rsh$ imply $r \nu \in \Rsh^+$. 

\np
In a previous paper \cite{DR10} we came across the analogue of the problem mentioned above: 
``Given a subset $\Ssh \subset \Rsh$ determine whether there is a positive system $\Rsh^+$ containing $\Ssh$''. 
An obvious necessary 
and sufficient condition (equivalent to the existence of the linear function $\varphi$ above) 
for the existence of a positive system $\Rsh^+$ containing $\Ssh$ is the requirement 
that the semigroup generated by $\Ssh$ does not contain $0$. Unfortunately, this combinatorial condition is 
not easy to verify. 
On the other hand, in our intended application in \cite{DR10},  the  condition $(\Sym^\cdot(\Msh))^{\gs} = \CC$ where 
$\Msh = \oplus_{\nu \in \Ssh}\, \gg^\nu$, arose naturally in the context of Geometric Invariant Theory.
This latter condition is necessary for the existence of a positive system $\Rsh^+$ as above.
To see this, note that $(\Sym^{\cdot}(\Msh))^{\gs}$ always contains at least the constants $\CC$, 
the inclusion $\gt \subset \gs$ implies $(\Sym^\cdot(\Msh))^{\gs} \subset (\Sym^\cdot(\Msh))^{\gt}$, and 
the condition that the semigroup generated by $\Ssh$ does not contain $0$ is
equivalent to $(\Sym^\cdot(\Msh))^{\gt} = \CC$.

\np
In fact, there is a stronger necessary condition for $\Ssh$ to be contained in a positive system. Since $\Rsh^+$ is saturated, if $\Ssh \subset \Rsh^+$
then $\overline{\Ssh} \subset \Rsh^+$, where $\overline{\Ssh}$ denotes the saturation of $S$, 
i.e., $\overline{\Ssh} = \QQ_+ \Ssh \cap \Rsh$.
Set $\overline{\Msh} := \oplus_{\nu \in \overline{\Ssh}}\, \gg^\nu$. It is easy to see that 
$(\Sym^\cdot(\Msh))^{\gt} = \CC$ if and only if $(\Sym^\cdot(\overline{\Msh}))^{\gt} = \CC$ 
and that we have the inclusions 
$(\Sym^\cdot(\Msh))^{\gs} \subset (\Sym^\cdot(\overline{\Msh}))^{\gs} \subset (\Sym^\cdot(\overline{\Msh}))^{\gt}$. 
In other words, if $\Ssh$ is contained in a positive system then $(\Sym^\cdot(\overline{\Msh}))^{\gs} = \CC$.

\np The goal of this note is to investigate whether either of the conditions $(\Sym^\cdot(\Msh))^{\gs} = \CC$ or 
$(\Sym^\cdot(\overline{\Msh}))^{\gs} = \CC$ is sufficient for the existence of a positive system $\Rsh^+$ containing $\Msh$.
It turns out that $(\Sym^\cdot(\Msh))^{\gs} = \CC$ is sufficient if and only if $\gg$ is of type $\A$ or $\D$ and 
$(\Sym^\cdot(\overline{\Msh}))^{\gs} = \CC$ is sufficient if and only if $\gg$ is classical or $\gg = \G_2$. 

\np
Using the connection between positive systems and linear functions $\varphi$ (valid when $\gt=\Zsh(\gs)$),
finding a positive system containing $\Msh$ 
is the same as finding a parabolic subalgebra $\gp_{\Msh}$ containing $\Msh$ with reductive part $\gs$, and
we will state our main result in this form.  
We will also state whether $\Ssh$ is saturated or not, rather than using the notation $\overline{\Msh}$.  
In the general case when $\gt\neq \Zsh(\gs)$, the existence of positive systems containing $\Msh$ is not equivalent 
to the existence of such a parabolic $\gp_{\Msh}$.  
However, our result, as stated below in terms of $\gp_{\Msh}$, is still valid in this case.

\np
\point
{\bf Main Theorem:} 
Let $\gg$ be a simple Lie algebra, $\gt\subset\gg$ a toral subalgebra, $\gs$ the centralizer of $\gt$,
$\Rsh$ the set of $\gt$-roots, $\Ssh\subset\Rsh$, and set $\Msh = \oplus_{\nu \in \Ssh}\, \gg^\nu$.   

\AlphaList
\begin{enumerate}
\item 
Assume that $(\Sym^{\cdot}(\Msh))^{\gs} = \CC$. 
If $\gg$ is of type $\A$ or $\D$ or if $\Ssh$ is saturated
and $\gg$ is of type $\B$, $\C$, or $\G_2$ then
there exists a parabolic subalgebra $\gp_{\Msh}$ with reductive part $\gs$ such that $\Msh\subset \gp_{\Msh}$.
\item If $\gg$ is not of type $\A$ or $\D$, there exist $\Ssh$ 
satisfying the condition that $(\Sym^\cdot(\Msh))^{\gs} = \CC$ such that no such parabolic $\gp_{\Msh}$ exists.
Moreover, if $\gg$ is $\F_{\!4}$, $\E_6$, $\E_7$, or $\E_8$, then $\Ssh$ can be chosen to be saturated. 
\end{enumerate}
\NumberList

\np
\point {\bf Reduction to $\bm{\gt=\Zsh(\gs)}.$}
In the main theorem we do not require that $\gt = \Zsh(\gs)$. However, the general case reduces to 
this case as follows: Set $\gt' := \Zsh(\gs)$ and let $\Rsh'$ be the set of $\gt'$-roots.  
The natural projection $\pi : (\gt')^* \to \gt^*$ 
induces a surjection of $\Rsh'$ onto $\Rsh$. Set $\Ssh':= \pi^{-1} (\Ssh)$ and notice that 
\[\Msh = \oplus_{\nu \in \Ssh}\, \gg^\nu = \oplus_{\nu' \in \Ssh'}\, \gg^{\nu'},\]
and that if $\Ssh$ is saturated, so is $\Ssh'$.  Moreover, the centralizer of $\gt'$ is again $\gs$.
Thus in proving that $(\Sym^{\cdot}(\Msh))^{\gs} = \CC$ is a sufficient condition we may assume that $\gt = \Zsh(\gs)$.
In the cases when we are proving that $(\Sym^{\cdot}(\Msh))^{\gs} = \CC$ is not sufficient, 
we provide examples in which $\gt = \Zsh(\gs)$.

\np 
For the rest of the paper we assume that $\gt = \Zsh(\gs)$.

\np 
\point
{\bf Organization and Conventions.} 
%The paper is organized as follows. 
In section \ref{sec:roots-and-spaces-classical} we describe explicitly all $\gt$-root systems and the 
respective $\gt$-root spaces for each of the classical simple Lie algebras. 
In section \ref{sec:proof-classical} we first prove the existence of $\gp_{\Msh}$ when $\gg$ is classical and $\Ssh$ is 
saturated.
We then handle the case of non-saturated $\Ssh$ in types $\A$ and $\D$, and finish the section by
giving examples in types $\B$ and $\C$ of non-saturated $\Ssh$ satisfying the condition
$(\Sym^{\cdot}(\Msh))^{\gs}=\CC$ for which no parabolic subalgebra $\gp_{\Msh}$ exists. 
In section \ref{sec:proof-exceptional} we first treat the case when $\gg$ is of type $\G_2$, proving the
result when $\Ssh$ is saturated and giving an example where $\Ssh$ is non-saturated.  We then construct
examples in types $\F_{\!4}$, $\E_6$, $\E_7$, and $\E_8$  of saturated $\Ssh$ for which 
$(\Sym^{\cdot}(\Msh))^{\gs}=\CC$ and for which no parabolic subalgebra $\gp_{\Msh}$ exists.
That is, in section \ref{sec:proof-classical} we establish all parts of the theorem dealing with classical
Lie algebras, and in section \ref{sec:proof-exceptional} we establish all parts dealing with the exceptional
Lie algebras. 

\np
Throughout the paper we work over the field of complex numbers $\CC$. 
All Lie algebras, modules, etc., are over $\CC$ unless explicitly stated otherwise.
The notation $\subset$ includes the possibility of equality.

\section{$\gt$-roots and $\gt$-root spaces for classical Lie algebras $\gg$.}
\label{sec:roots-and-spaces-classical}

\np
\point
First we describe the parabolic subalgebras and the corresponding sets $\Rsh$ for the classical Lie algebras. 
For convenience of notation we will work with the reductive Lie algebra $\ggl_n$ instead of $\gsl_n$. 
For the rest of this section
$\gg$ is a classical simple Lie algebra of type $\B$, $\C$, or $\D$ or $\gg = \ggl_n$. 
Moreover, we fix a Cartan subalgebra $\gh \subset \gg$. For a comprehensive source on simple complex Lie algebras 
we refer the reader to \cite{Bo}. For a treatment of parabolic subalgebras of $\gg$ containing a 
fixed Cartan subalgebra $\gh$, the reader may also consult \cite{DP}.

\np
\point
Let $\cP = \{\I_1, \ldots, \I_k\}$ be a partition of $\{1, \ldots, n\}$.   We say that $\cP$ is 
{\it totally ordered} if we have given a total order on the set $\{\I_1,\ldots, \I_k\}$.
We write $\cP(i)$ for the part of $\cP$ which contains $i$. The inequalities
$\cP(i) \prec \cP(j)$ and $\cP(i) \preceq \cP(j)$ are taken in the total order of the parts of $\cP$.
For the standard basis $\{ \vep_1, \ldots, \vep_n\}$ of $\gh^*$ we denote the dual basis of $\gh$ by
$\{h_1, \ldots, h_n\}$. 
A total order on the set $\{ \pm \delta_1, \ldots, \pm \delta_k\}$ is {\it compatible with multiplication by $-1$} if, 
for $x,y \in \{ \pm \delta_1, \ldots, \pm \delta_k\}$, $x \prec y$ implies $-y \prec -x$.
To simplify notation we adopt the convention that
$\B_1$, respectively $\C_1$, is a subalgebra of $\gg = \B_n$, respectively $\gg = \C_n$, isomorphic to 
$\A_1$ and whose roots are short, respectively long roots, of $\gg$. The subalgebras $\D_2=\A_1\oplus\A_1$ and 
$\D_3=\A_3$ of $\D_n$ have similar meaning.  

\np
Let $\gg$ be of type $\X_{n}=\A_n$, $\B_n$, $\C_n$, or $\D_n$ and let $\gs$ be a subalgebra of $\gg$ which 
is the reductive part of a parabolic subalgebra of $\gg$.  Every simple ideal of $\gs$ is isomorphic to $\A_r$
or $\X_r$ for some $r$.  Furthermore, if $\gg$ is not of type $\A_n$, $\gs$ has at most 
one simple ideal of type $\X_r$. 
For $\gg$ of type $\X_n=\B_n$, $\C_n$, or $\D_n$ the parabolic subalgebras of $\gg$
are split into two types depending on whether their reductive parts contain (Type II) or do not contain (Type I)
a simple ideal of type $\X_r$ (including $\B_1$, $\C_1$, $\D_2$, or $\D_3$).

\np
In the description of the combinatorics of the simple classical Lie algebras below,
the formulas for their parabolic subalgebras $\gp$ containing a fixed reductive part $\gs$ look very uniform
(e.g. \Old{11}).
In some instances this is misleading since the formulas do not explicitly indicate the subalgebra $\gs$ which, 
however, is an integral part of the structure of $\gp$.

\np
%{\bf Roots and parabolic subsets of the classical Lie algebras.} 
We now list the combinatorial descriptions of the parabolic subalgebras and related data in the classical cases.

\noindent
\point {$\bm{\gg=\ggl_n}$\bf.}
%\underline{$\gg=\ggl_n$}

\begin{enumerate}
\item The roots of $\gg$ are: $ \Delta = \{ \vep_i - \vep_j \, | \, 1 \leq i \neq j \leq n \} $

\item
Parabolic subalgebras of $\gg$ are in one-to-one correspondence with: \\

\begin{centering}
totally ordered partitions $\cP = (\I_1, \ldots, \I_k)$ of $\{1,\ldots,n\}$. \\
\end{centering}
\end{enumerate}

\np
Given a totally ordered partition $\cP$,

\begin{enumerate}
\setcounter{enumi}{2}
\item The roots of $\gp_\cP$ are $\{ \vep_i - \vep_j \, | \, i \neq j, \cP(i) \preceq \cP(j) \}$  
\item The roots of $\gs_\cP$ are $\{ \vep_i - \vep_j \, | \, i \neq j, \cP(i) = \cP(j) \}$
\item $\gs_\cP = \oplus_i \, \gs_\cP^i$, where $\gs_\cP^i \cong \ggl_{|\I_i|}$; 
\item The Cartan subalgebra of $\gs_\cP^i$ is spanned by $\{h_j\}_{j \in \I_i}$  
\item The roots of $\gs_\cP^i$ are $\{\vep_j - \vep_l \, | \, j \neq l \in \I_i\}$.   
\item $\gt_\cP$ has a basis $\{t_1, \ldots, t_k\}$ with $t_i = \frac{1}{|\I_i|} \sum_{j \in \I_i} h_j$ 
\end{enumerate}

If $\{\delta_1, \ldots, \delta_k\}$ is the basis of $\gt^*$ dual to $\{t_1, \ldots, t_k\}$ then

\begin{enumerate}
\setcounter{enumi}{8}
\item $\Rsh = \{ \delta_i - \delta_j \, | \, 1 \leq i \neq j \leq k\}.$
\item For $\nu = \delta_i - \delta_j \in \Rsh$, $\gs_\cP$-module $\gg^\nu\cong \V_i \otimes \V_j^*$, 
where $\V_i$ and $\V_j^*$ are the natural $\gs_\cP^i$-module and the dual of the natural $\gs_\cP^j$-module respectively,  all other factors of $\gs_{\cP}$ acting trivially.

\medskip
\item
The parabolic subalgebras of $\gg$ whose reductive part is $\gs_\cP$ are in a bijection with the ordered partitions $\cQ$ of
$\{1, \ldots, n\}$ whose parts are the same as the parts of $\cP$ or, equivalently, with total orders on the set 
$\{\delta_1, \ldots,
\delta_k\}$. 
\end{enumerate}

\bigskip
\noindent
\point {$\bm{\gg=\B_{n}}$}
%\underline{$\gg=\B_{n}$}

\begin{enumerate}
\item The roots of $\gg$ are: 
$\Delta = \{ \pm \vep_i \pm \vep_j, \pm \vep_i  \, | \, 1 \leq i \neq j \leq n \}.$

\item
Parabolic subalgebras of $\gg$ are in one-to-one correspondence with: \\

\begin{centering}
{\bf Type I:} 
\begin{minipage}[t]{0.75\textwidth}
pairs $(\cP, \sigma)$, where $\cP = (\I_1, \ldots, \I_k)$ is a totally ordered partition 
of $\{1, \ldots, n\}$ and $\sigma\colon \{1, \ldots, n \} \to \{\pm 1\}$ is a choice of signs. 
\end{minipage} \\

\bigskip
{\bf Type II:} 
\begin{minipage}[t]{0.75\textwidth}
pairs $(\cP, \sigma)$, where $\cP = (\I_0, \I_1, \ldots, \I_k)$ is a totally ordered partition 
of $\{1, \ldots, n\}$ with largest element $\I_0$ and 
$\sigma\colon \{1, \ldots, n \} \backslash \I_0 \to \{\pm 1\}$ is a choice of signs. 
\end{minipage}\\
\end{centering}
\end{enumerate}

\np
\und{In T}yp\und{e I:}

\begin{enumerate}
\setcounter{enumi}{2}
\item The roots of $\gp_{(\cP,\sigma)}$ are
$$\left\{{ \sigma(i) \vep_i - \sigma(j) \vep_j \, \st \, i \neq j, \cP(i) \preceq \cP(j) }\right\} \cup \{\sigma(i) \vep_i + \sigma(j) \vep_j, \sigma(i) \vep_i \, \st \, i \neq j\}$$
\item The roots of $\gs_{(\cP,\sigma)}$ are $\{ \sigma(i)\vep_i - \sigma(j)\vep_j \, | \, i \neq j, \cP(i) = \cP(j) \}$.
\item $\gs_{(\cP,\sigma)} = \oplus_i \, \gs_{(\cP,\sigma)}^i$, where $\gs_{(\cP,\sigma)}^i \cong \ggl_{|\I_i|}$.
\item The Cartan subalgebra of $\gs_\cP^i$ is spanned by $\{\sigma(j) h_j\}_{j \in \I_i}$ 
\item The roots of $\gs_\cP^i$ are $\{\sigma(j)\vep_j - \sigma(l)\vep_l \, | \, j \neq l \in \I_i\}$.  
\item 
$\gt_{(\cP,\sigma)}$ has a basis $\{t_1, \ldots, t_k\}$ with $t_i = \frac{1}{|\I_i|} \sum_{j \in \I_i} \sigma(j) h_j$.
\end{enumerate}

If $\{\delta_1, \ldots, \delta_k\}$ the basis of $\gt^*$ dual to $\{t_1, \ldots, t_k\}$ then

\begin{enumerate}
\setcounter{enumi}{8}
\item 
$\Rsh = \{ \pm \delta_i \pm \delta_j, \pm \delta_i \, | \, 1 \leq i \neq j \leq k\} \cup \{\pm 2 \delta_i \, | \, |\I_i|>1\}.$
\item For $\nu \in \Rsh$, 

\begin{enumerate}
\item $\gg^\nu \cong \V_i^\pm \otimes \V_j^\pm$ if $\nu = \pm \delta_i \pm \delta_j$,
\item $\gg^\nu\cong\V_i^\pm$ if $\nu = \pm \delta_i$, and 
\item $\gg^\nu\cong\Lambda^2 \V_i^\pm$ if $\nu = \pm 2 \delta_i$, 
\end{enumerate}

\np
where $\V_i^+$ and $\V_i^-$ respectively are the natural $\gs_{(\cP,\sigma)}^i$-module and its dual, and
all other factors of $\gs_{(\cP,\sigma)}$ act trivially.

\medskip
\item
The parabolic subalgebras of $\gg$ whose reductive part is $\gs_{\cP,\sigma}$ are in a bijection with the pairs $(\cQ, \tau)$ 
such that the parts of $\cQ$ are the same as the parts of $\cP$ and $\sigma_{|\I_i} = \pm \tau_{|\I_i}$ for every part $\I_i$ or, 
equivalently, with total orders on the set $\{\pm \delta_1, \ldots, \pm \delta_k\}$ compatible with multiplication by $-1$.
\end{enumerate}

\np
\und{In T}yp\und{e II:}

\begin{enumerate}
\setcounter{enumi}{2}
\item The roots of $\gp_{(\cP,\sigma)}$ are 

\medskip
$
\{ \sigma(i) \vep_i - \sigma(j) \vep_j \, | \,  i \neq j, \cP(i) \preceq \cP(j) \prec \I_0\} 
\cup \{\pm \vep_i \pm \vep_j, \pm \vep_i \, | \, i \neq j , i \in \I_0, j \in \I_0\} $ \hfill

\medskip
\hfill $\cup
\{\sigma(i) \vep_i + \sigma(j) \vep_j, \sigma(i) \vep_i \, | \, i \neq j, i \not \in \I_0,j \not \in \I_0\} \cup 
\{\sigma(i) \vep_i \pm \vep_j \, | \, i \not \in \I_0, j \in \I_0\}$ 

\item The roots of $\gs_{(\cP,\sigma)}$ are
$$\{ \sigma(i)\vep_i - \sigma(j)\vep_j \, | \, i \neq j, \cP(i) = \cP(j) \prec \I_0 \} \cup \{\pm \vep_i \pm \vep_j, \pm \vep_i \, | \, i \neq j \in \I_0\}.$$
\item $\gs_{(\cP,\sigma)} = \oplus_i \, \gs_{(\cP,\sigma)}^i$, 
where $\gs_{(\cP,\sigma)}^0 \cong \B_{|\I_0|}$ and $\gs_{(\cP,\sigma)}^i \cong \ggl_{|\I_i|}$ for $i>0$.
\item The Cartan subalgebra of $\gs_\cP^i$ is spanned by $\{h_j\}_{j \in \I_0}$ for $i=0$ and 
$\{\sigma(j) h_j\}_{j \in \I_i}$ for $i>0$. 
\item The roots of $\gs_\cP^i$ are $\{\pm \vep_j \pm \vep_l, \pm \vep_j \, | \, j \neq l \in \I_0\}$ for $i = 0$ and
$\{\sigma(j)\vep_j - \sigma(l)\vep_l \, | \, j \neq l \in \I_i\}$ for $i>0$.  
\item 
$\gt_{(\cP,\sigma)}$ has a basis $\{t_1, \ldots, t_k\}$ with $t_i = \frac{1}{|\I_i|} \sum_{j \in \I_i} \sigma(j) h_j$.
\end{enumerate}

\np
If $\{\delta_1, \ldots, \delta_k\}$ is the basis of $\gt^*$ dual to $\{t_1, \ldots, t_k\}$ then

\begin{enumerate}
\setcounter{enumi}{8}
\item 
$\Rsh = \{ \pm \delta_i \pm \delta_j, \pm \delta_i \, | \, 1 \leq i \neq j \leq k\} \cup \{\pm 2 \delta_i \, | \, |\I_i|>1\}.$
\item For $\nu \in \Rsh$, 

\begin{enumerate}
\item $\gg^\nu\cong \V_i^\pm \otimes \V_j^\pm$ if $\nu = \pm \delta_i \pm \delta_j$, 
\item $\gg^\nu\cong \V_i^\pm \otimes \V_0$ if $\nu = \pm \delta_i$,  and 
\item $\gg^\nu\cong\Lambda^2 \V_i^\pm$ if $\nu = \pm 2 \delta_i$ 
\end{enumerate}

\np
where
$\V_i^+$ and $\V_i^-$ denote the natural $\gs_{(\cP,\sigma)}^i$-module and its dual respectively for $i>0$,
$\V_0$ denotes the natural $\gs_{(\cP,\sigma)}^0$-module,
and all other factors of $\gs_{(\cP,\sigma)}$ act trivially. Note that, if $\gs_{(\cP,\sigma)} = \B_1 \cong \gsl_2$, then $\V_0$ is the three dimensional
irreducible $\gs_{(\cP,\sigma)}$-module.

\medskip
\item
The parabolic subalgebras of $\gg$ whose reductive part is $\gs_{\cP,\sigma}$ are in a bijection with the pairs $(\cQ, \tau)$ 
such that the parts of $\cQ$ are the same as the parts of $\cP$, $\I_0$ is the largest element of $\cQ$,
and $\sigma_{|\I_i} = \pm \tau_{|\I_i}$ for every part $\I_i \neq \I_0$ or, 
equivalently, with total orders on the set $\{\pm \delta_1, \ldots, \pm \delta_k\}$ compatible with multiplication by $-1$. 
\end{enumerate}

\bigskip
\noindent
\point {$\bm{\gg=\C_{n}}$}
%\underline{$\gg=\C_{n}$}

\begin{enumerate}
\item The roots of $\gg$ are: 
$\Delta = \{ \pm \vep_i \pm \vep_j, \pm 2 \vep_i  \, | \, 1 \leq i \neq j \leq n \}.$

\item
Parabolic subalgebras of $\gg$ are in one-to-one correspondence with: \\

\begin{centering}
{\bf Type I:} 
\begin{minipage}[t]{0.75\textwidth}
pairs $(\cP, \sigma)$, where $\cP = (\I_1, \ldots, \I_k)$ is a totally ordered partition 
of $\{1, \ldots, n\}$ and $\sigma\colon \{1, \ldots, n \} \to \{\pm 1\}$ is a choice of signs. 
\end{minipage} \\

\bigskip
{\bf Type II:} 
\begin{minipage}[t]{0.75\textwidth}
pairs $(\cP, \sigma)$, where $\cP = (\I_0, \I_1, \ldots, \I_k)$ is a totally ordered partition 
of $\{1, \ldots, n\}$ with largest element $\I_0$ and 
$\sigma\colon \{1, \ldots, n \} \backslash \I_0 \to \{\pm 1\}$ is a choice of signs. 
\end{minipage}\\
\end{centering}
\end{enumerate}

\np
\und{In T}yp\und{e I:}

\begin{enumerate}
\setcounter{enumi}{2}
\item The roots of $\gp_{(\cP,\sigma)}$  are
$$\{ \sigma(i) \vep_i - \sigma(j) \vep_j \, | \, i \neq j, \cP(i) \preceq \cP(j) \} \cup \{\sigma(i) \vep_i + \sigma(j) \vep_j, 2 \sigma(i) \vep_i \, | \, i \neq j\}$$
\item The roots of $\gs_{(\cP,\sigma)}$ are $\{ \sigma(i)\vep_i - \sigma(j)\vep_j \, | \, i \neq j, \cP(i) = \cP(j) \}$.
\item $\gs_{(\cP,\sigma)} = \oplus_i \, \gs_{(\cP,\sigma)}^i$, where $\gs_{(\cP,\sigma)}^i \cong \ggl_{|\I_i|}$.
\item The Cartan subalgebra of $\gs_\cP^i$ is spanned by $\{\sigma(j) h_j\}_{j \in \I_i}$.
\item The roots of $\gs_{(\cP,\sigma)}^i$ are $\{\sigma(j)\vep_j - \sigma(l)\vep_l \, | \, j \neq l \in \I_i\}$.  
\item 
$\gt_{(\cP,\sigma)}$ has a basis $\{t_1, \ldots, t_k\}$ with $t_i = \frac{1}{|\I_i|} \sum_{j \in \I_i} \sigma(j) h_j$.
\end{enumerate}

\noindent
If $\{\delta_1, \ldots, \delta_k\}$ is the basis of $\gt^*$ dual to $\{t_1, \ldots, t_k\}$ then

\begin{enumerate}
\setcounter{enumi}{8}
\item 
$\Rsh = \{ \pm \delta_i \pm \delta_j, \pm 2 \delta_i \, | \, 1 \leq i \neq j \leq k\}.$
\item For $\nu \in \Rsh$, 

\begin{enumerate}
\item $\gg^\nu\cong \V_i^\pm \otimes \V_j^\pm$ if  $\nu = \pm \delta_i \pm \delta_j$.
\item $\gg^\nu\cong\Sym^2 \V_i^\pm$ if for $\nu = \pm 2 \delta_i$.
\end{enumerate}

\np
where $\V_i^+$ and $\V_i^-$ are the natural $\gs_{(\cP,\sigma)}^i$-module and its dual, and
all other factors of $\gs_{(\cP,\sigma)}$ act trivially.

\medskip
\item
The parabolic subalgebras of $\gg$ whose reductive part is $\gs_{\cP,\sigma}$ are in a bijection with the pairs $(\cQ, \tau)$ 
such that the parts of $\cQ$ are the same as the parts of $\cP$ and $\sigma_{|\I_i} = \pm \tau_{|\I_i}$ for every part $\I_i$ or, 
equivalently, with total orders on the set $\{\pm \delta_1, \ldots, \pm \delta_k\}$ compatible with multiplication by $-1$.
\end{enumerate}

\np
\und{In T}yp\und{e II:}

\begin{enumerate}
\setcounter{enumi}{2}
\item The roots of $\gp_{(\cP,\sigma)}$ are

\medskip
$\{ \sigma(i) \vep_i - \sigma(j) \vep_j \, | \,  i \neq j, \cP(i) \preceq \cP(j) \prec \I_0\} 
\cup \{\pm \vep_i \pm \vep_j, \pm 2 \vep_i \, | \, i \neq j, i \in \I_0, j \in \I_0\} \cup$ 

\medskip
\hfill $
\{\sigma(i) \vep_i + \sigma(j) \vep_j, \sigma(i) 2 \vep_i \, | \, i \neq j, i \not \in \I_0, j \not \in \I_0\} \cup 
\{\sigma(i) \vep_i \pm \vep_j \, | \, i \not \in \I_0, j \in \I_0\}$ 

\medskip
\item The roots of $\gs_{(\cP,\sigma)}$ are

\medskip
\begin{centering}
$\{ \sigma(i)\vep_i - \sigma(j)\vep_j \, | \, i \neq j, \cP(i) = \cP(j) \prec \I_0 \} \cup \{\pm \vep_i \pm \vep_j, \pm 2 \vep_i \, | \, i \neq j \in \I_0\}.$ \\
\end{centering}
\medskip
\item 
$\gs_{(\cP,\sigma)} = \oplus_{i=0}^k \, \gs_{(\cP,\sigma)}^i$, where $\gs_{(\cP,\sigma)}^0 \cong \C_{|\I_0|}$ 
and $\gs_{(\cP,\sigma)}^i \cong \ggl_{|\I_i|}$ for $i>0$.
\item The Cartan subalgebra of $\gs_{(\cP,\sigma)}^i$ is spanned by $\{h_j\}_{j \in \I_0}$ for $i=0$ and 
$\{\sigma(j) h_j\}_{j \in \I_i}$ for $i>0$. 
\item The roots of $\gs_{(\cP,\sigma)}^i$ are
$\{\pm \vep_j \pm \vep_l, \pm 2 \vep_j \, | \, j \neq l \in \I_0\}$ for $i = 0$ and
$\{\sigma(j)\vep_j - \sigma(l)\vep_l \, | \, j \neq l \in \I_i\}$ for $i>0$.  
\item 
$\gt_{(\cP,\sigma)}$ has a basis $\{t_1, \ldots, t_k\}$ with $t_i = \frac{1}{|\I_i|} \sum_{j \in \I_i} \sigma(j) h_j$ 
\end{enumerate}

\noindent
If $\{\delta_1, \ldots, \delta_k\}$ is the basis of $\gt^*$ dual to $\{t_1, \ldots, t_k\}$ then

\begin{enumerate}
\setcounter{enumi}{8}
\item $\Rsh = \{ \pm \delta_i \pm \delta_j, \pm \delta_i, \pm 2 \delta_i \, | \, 1 \leq i \neq j \leq k\}$.
\item For $\nu \in \Rsh$, 

\begin{enumerate}
\item $\gg^\nu\cong\V_i^\pm \otimes \V_j^\pm$ if $\nu = \pm \delta_i \pm \delta_j$,
\item $\gg^\nu\cong\V_i^\pm \otimes \V_0$ if $\nu = \pm \delta_i$,
\item $\gg^\nu$ is isomorphic to $\Sym^2 \V_i^\pm$ if $\nu = \pm 2 \delta_i$,
\end{enumerate}

\np
where $\V_i^+$ and $\V_i^-$ denote the natural $\gs_{(\cP,\sigma)}^i$-module and its dual 
for $i>0$, 
$\V_0$ is the natural $\gs_{(\cP,\sigma)}^0$-module, and where 
all other factors of $\gs_{(\cP,\sigma)}$ act trivially.
Note that, if $\gs_{(\cP,\sigma)} = \C_1 \cong \gsl_2$, then $\V_0$ is the two dimensional
irreducible $\gs_{(\cP,\sigma)}$-module.

\medskip
\item
The parabolic subalgebras of $\gg$ whose reductive part is $\gs_{\cP,\sigma}$ are in a bijection 
with the pairs $(\cQ, \tau)$ such that the parts of $\cQ$ are the same as the parts 
of $\cP$, $\I_0$ is the largest element of $\cQ$, 
and $\sigma_{|\I_i} = \pm \tau_{|\I_i}$ for every part $\I_i \neq \I_0$ or, equivalently, with total orders on 
the set $\{\pm \delta_1, \ldots, \pm \delta_k\}$ compatible with multiplication by $-1$. 
\end{enumerate}

\bigskip
\noindent
\point {$\bm{\gg=\D_{n}}$\bf.}
%\underline{$\gg=\D_{n}$}

\begin{enumerate}
\item The roots of $\gg$ are: 
$\Delta = \{ \pm \vep_i \pm \vep_j,  \, | \, 1 \leq i \neq j \leq n \}.$

\item
Parabolic subalgebras of $\gg$ are determined by: \\

\begin{centering}
{\bf Type I:} 
\begin{minipage}[t]{0.75\textwidth}
pairs $(\cP, \sigma)$, where $\cP = (\I_0, \I_1, \ldots, \I_k)$ is a totally ordered partition 
of $\{1, \ldots, n\}$ and $\sigma: \{1, \ldots, n \} \to \{\pm 1\}$ is a choice of signs. 

\np
Two pairs $(\cP', \sigma')$ and $(\cP'', \sigma'')$
determine the same parabolic subalgebra if and only if $\cP'$ and $\cP''$ are the same ordered 
partitions whose maximal part $\I_0$ contains one element and $\sigma'$ and $\sigma''$ 
coincide on $\{1, \ldots, n\} \backslash \I_0$.
\end{minipage} \\

\bigskip
{\bf Type II:} 
\begin{minipage}[t]{0.75\textwidth}
pairs $(\cP, \sigma)$, where $\cP = (\I_0, \I_1, \ldots, \I_k)$ is a totally ordered partition 
of $\{1, \ldots, n\}$ with largest element $\I_0$ such that $|\I_0| \geq 2$ and 
$\sigma: \{1, \ldots, n \} \backslash \I_0 \to \{\pm 1\}$ is a choice of signs. 
\end{minipage}\\
\end{centering}
\end{enumerate}

\np
\und{In T}yp\und{e I:}

\begin{enumerate}
\setcounter{enumi}{2}
\item The roots of $\gp_{(\cP,\sigma)}$ are  
$\{ \sigma(i) \vep_i - \sigma(j) \vep_j \, | \, i \neq j, \cP(i) \preceq \cP(j) \} \cup \{\sigma(i) \vep_i + \sigma(j) \vep_j \, | \, i \neq j\}$ 
\item The roots of $\gs_{(\cP,\sigma)}$ are $\{ \sigma(i)\vep_i - \sigma(j)\vep_j \, | \, i \neq j, \cP(i) = \cP(j) \}$.
\item $\gs_{(\cP,\sigma)} = \oplus_i \, \gs_{(\cP,\sigma)}^i$, where $\gs_{(\cP,\sigma)}^i \cong \ggl_{|\I_i|}$.
\item The Cartan subalgebra of $\gs_\cP^i$ is spanned by $\{\sigma(j) h_j\}_{j \in \I_i}$ 
\item The roots of $\gs_{(\cP,\sigma)}^i$ are $\{\sigma(j)\vep_j - \sigma(l)\vep_l \, | \, j \neq l \in \I_i\}$.  
\item 
$\gt_{(\cP,\sigma)}$ has a basis $\{t_1, \ldots, t_k\}$ with $t_i = \frac{1}{|\I_i|} \sum_{j \in \I_i} \sigma(j) h_j$.
\end{enumerate}

\noindent
If $\{\delta_1, \ldots, \delta_k\}$ is the basis of $\gt^*$ dual to $\{t_1, \ldots, t_k\}$ then

\begin{enumerate}
\setcounter{enumi}{8}
\item 
$\Rsh = \{ \pm \delta_i \pm \delta_j | \, 1 \leq i \neq j \leq k\} \cup \{ \pm 2 \delta_i \, | \, |\I_i| >1\}.$
\item For $\nu \in \Rsh$, 

\begin{enumerate}
\item $\gg^\nu\cong \V_i^\pm \otimes \V_j^\pm$ if  $\nu = \pm \delta_i \pm \delta_j$.
\item $\gg^\nu\cong\Lambda^2 \V_i^\pm$ if $\nu = \pm 2 \delta_i$.
\end{enumerate}

\np
where $\V_i^+$ and $\V_i^-$ are the natural $\gs_{(\cP,\sigma)}^i$-module and its dual, and
all other factors of $\gs_{(\cP,\sigma)}$ act trivially.

\medskip
\item
Every parabolic subalgebra of $\gg$ whose reductive part is $\gs_{\cP,\sigma}$ corresponds to a pair $(\cQ, \tau)$ 
such that the parts of $\cQ$ are the same as the parts of $\cP$ and $\sigma_{|\I_i} = \pm \tau_{|\I_i}$ for every part $\I_i$ or, 
equivalently, to a total order on the set $\{\pm \delta_1, \ldots, \pm \delta_k\}$ compatible with multiplication by $-1$. 
Note that this correspondence is not bijective since
two different total orders may determine the same parabolic subalgebra.
\end{enumerate}

\np
\und{In T}yp\und{e II:}

\begin{enumerate}
\setcounter{enumi}{2}
\item Roots of $\gp_{(\cP,\sigma)} = $

\medskip
$\left\{ \sigma(i) \vep_i - \sigma(j) \vep_j \, \st  \,  i \neq j, \cP(i) \preceq \cP(j) \prec \I_0\right\} 
\cup \{\pm \vep_i \pm \vep_j \, | \, i \neq j , i \in \I_0, j \in \I_0\} $

\medskip
\hfill
$\cup
\{\sigma(i) \vep_i + \sigma(j) \vep_j \, | \, i \neq j, i \not \in \I_0, j \not \in \I_0\} \cup 
\{\sigma(i) \vep_i \pm \vep_j \, | \, i \not \in \I_0, j \in \I_0\}$ 

\medskip
\item Roots of $\gs_{(\cP,\sigma)}= 
\{ \sigma(i)\vep_i - \sigma(j)\vep_j \, | \, i \neq j, \cP(i) = \cP(j) \prec \I_0 \} \cup \{\pm \vep_i \pm \vep_j \, | \, i \neq j \in \I_0\}.$ 

\medskip
\item $\gs_{(\cP,\sigma)} = \oplus_{i=0}^k \, \gs_{(\cP,\sigma)}^i$, where $\gs_{(\cP,\sigma)}^0 \cong \D_{|\I_0|}$ and 
$\gs_{(\cP,\sigma)}^i \cong \ggl_{|\I_i|}$ for $i>0$.
\item 
Cartan subalgebra of $\gs_{(\cP,\sigma)}^i$ is spanned by $\{h_j\}_{j \in \I_0}$ for $i=0$ and by 
$\{\sigma(j) h_j\}_{j \in \I_i}$ for $i>0$; 

\medskip
\item 
roots of $\gs_{(\cP,\sigma)}^i$ are
$\{\pm \vep_j \pm \vep_l \, | \, j \neq l \in \I_0\}$ for $i = 0$ and
$\{\sigma(j)\vep_j - \sigma(l)\vep_l \, | \, j \neq l \in \I_i\}$ for $i>0$.  

\medskip
\item 
$\gt_{(\cP,\sigma)}$ has a basis $\{t_1, \ldots, t_k\}$ with $t_i = \frac{1}{|\I_i|} \sum_{j \in \I_i} \sigma(j) h_j$.
\end{enumerate}

\noindent
If $\{\delta_1, \ldots, \delta_k\}$ is the basis of $\gt^*$ dual to $\{t_1, \ldots, t_k\}$ then

\begin{enumerate}
\setcounter{enumi}{8}
\item 
$\Rsh = \{ \pm \delta_i \pm \delta_j, \pm \delta_i \, | \, 1 \leq i \neq j \leq k\} \cup \{\pm 2 \delta_i \, | \, |\I_i| >1\}.$

\item For $\nu \in \Rsh$, 

\begin{enumerate}
\item $\gg^\nu\cong \V_i^\pm \otimes \V_j^\pm$ if $\nu = \pm \delta_i \pm \delta_j$,
\item $\gg^\nu\cong\V_i^\pm \otimes \V_0$ if $\nu = \pm \delta_i$, 
\item $\gg^\nu\cong\Lambda^2 \V_i^\pm$ if $\nu = \pm 2 \delta_i$,
\end{enumerate}

\np
where $\V_i^+$ and $\V_i^-$ denote the natural $\gs_{(\cP,\sigma)}^i$-module and its dual 
for $i>0$, 
$\V_0$ is the natural  $\gs_{(\cP,\sigma)}^0$-module, and where 
all other factors of $\gs_{(\cP,\sigma)}$ act trivially.
Note that, if $\gs_{(\cP,\sigma)} = \D_2 \cong \gsl_2 \oplus \gsl_2$, then $\V_0$ is the (external) tensor product of 
two two-dimensional
irreducible $\gsl_2$-modules; if $\gs_{(\cP,\sigma)} = \D_3 \cong \gsl_4$, 
then $\V_0$ the six dimensional irreducible $\gs_{(\cP,\sigma)}$-module
which is the second exterior power of the natural representation of $\gsl_4$.

\medskip
\item
The parabolic subalgebras of $\gg$ whose reductive part is $\gs_{(\cP,\sigma)}$ are in a bijection with the pairs $(\cQ, \tau)$ 
such that the parts of $\cQ$ are the same as the parts of $\cP$, $\I_0$ is the largest element of $\cQ$, and $\sigma_{|\I_i} = \pm \tau_{|\I_i}$ 
for every part $\I_i$ or, 
equivalently, with total orders on the set $\{\pm \delta_1, \ldots, \pm \delta_k\}$. 
\end{enumerate}

\section{Proof of the Main Theorem when $\gg$ is classical.}
\label{sec:proof-classical}

\np
\point 
{\bf Existence of $\gp_{\Msh}$ when $\Ssh$ is saturated.}
The idea is simple: using $\Ssh$ we define a binary relation $\prec$ on the set $\{\delta_1, \ldots, \delta_k\}$
(respectively on $\{\pm \delta_1, \ldots, \pm \delta_k\}$) and using the fact that $(\Sym^\cdot(\Msh))^\gs = \CC$ 
we prove that $\prec$ can be extended to a total order (respectively, to a total order compatible with multiplication by $-1$).
The proof follows the same logic in all cases but is least technical in the case when $\gg = \ggl_n$. 
For clarity of exposition
we present the proof for $\gg= \ggl_n$ first. 
Thoughout the proof the partition $\cP$ (and the choice of signs $\sigma$) 
are fixed and instead of $\gs_\cP$ (or $\gs_{(\cP, \sigma)}$) and $\gs_\cP^i$ (or $\gs_{(\cP, \sigma)}^i$) we write 
$\gs$ and $\gs^i$ respectively.

\np 
First we consider the case when $\gg = \ggl_n$. Define a binary relation $\prec$ on $\{\delta_1, \delta_2, \ldots, \delta_k\}$ by setting

\begin{equation} \label{eq2}
\delta_i \prec \delta_j \quad  {\text { if }} \quad \nu = \delta_i - \delta_j \in \Ssh.
\end{equation}
The existence of a parabolic subalgebra $\gp_\Msh$ with reductive part $\gs$ and containing $\Msh$ is equivalent
the existence of a total order on $\{\delta_1, \delta_2, \ldots, \delta_k\}$ which extends $\prec$.

\np
Note that $\prec$ can be extended to a total order on $\{\delta_1, \delta_2, \ldots, \delta_k\}$ if and only if 
there is no cycle 

\begin{equation} \label{cycle}
\delta_{i_1} \prec \delta_{i_2} \prec \cdots \prec \delta_{i_l} \prec \delta_{i_1}.
\end{equation}

\np 
Assume that $\prec$ cannot  be extended to a total order on  $\{\delta_1, \cdots, \delta_k\}$ and consider a cycle
\eqref{cycle} of minimal length. Then $\nu_1 = \delta_{i_1} - \delta_{i_2}, \nu_2 = \delta_{i_2} - \delta_{i_3}, \cdots, \nu_l = \delta_{i_l} - \delta_{i_1}$ 
is a sequence of distinct elements of $\Ssh$.  Hence $\gg^{\nu_1} \oplus \gg^{\nu_2} \oplus \cdots \oplus \gg^{\nu_l}$
is a submodule of $\Msh$ and $\Sym^l (\gg^{\nu_1} \oplus \gg^{\nu_2} \oplus \cdots \oplus \gg^{\nu_l})$ is a submodule of $\Sym^\cdot(\Msh)$
containing $\gg^{\nu_1} \otimes \gg^{\nu_2} \otimes \cdots \otimes \gg^{\nu_l}$. On the other hand,

\begin{equation} \label{eq12}
%\rule{0.5cm}{0cm}
\begin{array}{ccl}
\gg^{\nu_1} \otimes\gg^{\nu_2} \otimes \cdots \otimes \gg^{\nu_l} & \cong &
(\V_{i_1} \otimes \V_{i_2}^*) \otimes (\V_{i_2} \otimes \V_{i_3}^*) \otimes \cdots \otimes (\V_{i_l} \otimes \V_{i_1}^*) \\ & \cong &
(\V_{i_1} \otimes \V_{i_1}^*) \otimes (\V_{i_2} \otimes \V_{i_2}^*) \otimes \cdots \otimes (\V_{i_l} \otimes \V_{i_l}^*), \end{array}
\end{equation}
where the lower index of a module shows which component of $\gs$ acts non-trivially on it. 
Since, for every $1 \leq j \leq l$, $\V_{i_j} \otimes \V_{i_j}^*$ contains the trivial $\gs^{i_j}$-module, 
(\ref{eq12}) shows that $\gg^{\nu_1} \otimes\gg^{\nu_2} \otimes \cdots \otimes \gg^{\nu_l}$ contains the trivial $\gs$-module 
which contradicts the assumption that
$(\Sym^\cdot(\Msh))^\gs = \CC$. This contradiction shows that $\prec$ can be extended to a total order on the 
set $\{\delta_1, \cdots, \delta_k\}$, which completes the proof when $\gg = \ggl_n$.

\np 
Next we consider the case when $\gg \neq \ggl_n$, i.e., we assume that $\gg$ is a simple classical Lie algebra not of type $\A$. 
Define a binary relation $\prec$ on $\{\pm, \delta_1, \pm \delta_2, \cdots, \pm \delta_k\}$ by setting

\begin{equation} \label{eq22} 
\begin{array}{lll}
s_i \delta_i \prec s_j \delta_j, \text{   }  i \neq j & {\text { if }} &  \nu = s_i \delta_i - s_j \delta_j \in \Ssh \\
s_i \delta_i \prec - s_i \delta_i & {\text { if }} & \nu = \left\{
\begin{array}{ccc}
s_i \delta_i & {\text{when}} & \gg = \B_n {\text{  or  }} \gg=\D_n, {\text{ Type II}}\\
2 s_i \delta_i  & {\text{when}} & \gg = \C_n {\text{  or  }} \gg=\D_n, {\text{ Type I}}
\end{array} \right. \in \Ssh,
\end{array}
\end{equation}
where $s_i, s_j = \pm$. Note that $\prec$ is compatible with multiplication by $-1$. 

\np
The existence of a parabolic subalgebra $\gp_\Msh$ with reductive part $\gs$ and containing $\Msh$ is equivalent
the existence of a total order on $\{\pm \delta_1, \pm \delta_2, \cdots, \pm \delta_k\}$ compatible
with multiplication by $-1$ which extends $\prec$.

\np
Note that $\prec$ can be extended to a total order on $\{\pm \delta_1, \pm \delta_2, \cdots, \pm \delta_k\}$ compatible
with multiplication by $-1$ if and only if there is no cycle 

\begin{equation} \label{eq31}
s_1 \delta_{i_1} \prec s_2 \delta_{i_2} \prec \cdots \prec s_l \delta_{i_l} \prec s_1 \delta_{i_1}
\end{equation}

\np 
Assume that $\prec$ cannot  be extended to a total order on $\{\pm \delta_1, \cdots, \pm \delta_k\}$ compatible with multiplication by $-1$ and 
consider a cycle \eqref{eq31} of minimal length. It gives rise to 
a sequence $\nu_1, \cdots, \nu_l \in \Ssh$ induced from (\ref{eq22}). 
More precisely, 
$$
\nu_j = \left\{ \begin{array}{ccl} 
s_j \delta_{i_j} - s_{j+1} \delta_{i_{j+1}} & {\text{ if }} &  \delta_{i_j} \neq \delta_{i_{j+1}}\\
s_j \delta_{i_j} & {\text { if }} &  \delta_{i_j} = \delta_{i_{j+1}}, \gg = \B_n {\text{  or  }} \gg=\D_n, {\text{ Type II}}\\
2s_j \delta_{i_j} & {\text { if }} &  \delta_{i_j} = \delta_{i_{j+1}}, \gg = \C_n {\text{  or  }} \gg=\D_n, {\text{ Type I}},
\end{array} \right. 
$$
where $s_{l+1} = s_1$ and $\delta_{i_{l+1}} = \delta_{i_1}$. 

\np
The minimality of \eqref{eq31} implies that every element $\nu$ of $\Rsh$ appears at most twice in the sequence $\nu_1, \nu_2, \ldots, n_l$.
Moreover, if $\nu= \pm \delta_i$ or $\nu = \pm 2 \delta_i$, then $\nu$ appears at most once in this sequence.

\np
First we consider the case when  $\delta_{i_j} \neq \delta_{i_{j+1}}$ for every $j$. 
In this case $\nu_j = s_j \delta_{i_j} - s_{j+1} \delta_{i_{j+1}}$ for every $j$. 
Let $\lambda_1,  \ldots, \lambda_s$ be the elements of $\Rsh$ that appear once in the sequence $\nu_1, \nu_2, \ldots, \nu_l$ and let
$\mu_1, \ldots, \mu_t$ be those that appear twice. Clearly, $l = s + 2t$. 
Moreover $\gg^{\lambda_1} \oplus \cdots \oplus \gg^{\lambda_s} \oplus \gg^{\mu_1} \oplus \cdots \oplus \gg^{\mu_t}$
is a submodule of $\Msh$ and $\Sym^l (\gg^{\lambda_1} \oplus \cdots \oplus \gg^{\lambda_s} \oplus \gg^{\mu_1} \oplus \cdots \oplus \gg^{\mu_t})$ 
is a submodule of $\Sym^\cdot(\Msh)$
containing 

\begin{equation} \label{eq411}
\gg^{\lambda_1} \otimes \cdots \otimes \gg^{\lambda_s} \otimes \Sym^2\gg^{\mu_1} \otimes \cdots \otimes \Sym^2 \gg^{\mu_t}.
\end{equation} 
We will prove that the $\gs$-module \eqref{eq411} contains the trivial $\gs$-module which, as in the case when $\gg = \ggl_n$, will complete the proof.

\np
Indeed, if $\mu_{j'} = \nu_j$ then

\begin{equation} \label{eq412}
\begin{array}{cl}
&\Sym^2 \gg^{\mu_{j'}}  =  \Sym^2 \gg^{\nu_j} = \Sym^2 (V_{i_j}^{s_j} \otimes \V_{i_{j+1}}^{-s_{j+1}}) = \\
 & \Sym^2 \V_{i_j}^{s_j} \otimes \Sym^2 \V_{i_{j+1}}^{-s_{j+1}} \oplus  \Lambda^2 \V_{i_j}^{s_j} \otimes \Lambda^2 \V_{i_{j+1}}^{-s_{j+1}} \supset 
 \Sym^2 \V_{i_j}^{s_j} \otimes \Sym^2 \V_{i_{j+1}}^{-s_{j+1}}.
 \end{array}
\end{equation}
Replacing in \eqref{eq411} each term of the form $\Sym^2 \gg^{\mu_{j'}}$ with the corresponding term
$\Sym^2 \V_{i_j}^{s_j} \otimes \Sym^2 \V_{i_{j+1}}^{-s_{j+1}}$ from \eqref{eq412}, we obtain another submodule of \eqref{eq411}.
This latest submodule is a tensor product of factors of the form $\V_{i_j}^\pm$ and $\Sym^2 \V_{i_j}^\pm$. Moreover, the component $\V_i$ appears
in one of the following groups:
\[\V_i^+ \otimes \V_i^+ \otimes \V_i^- \otimes \V_i^-,\, 
\V_i^+ \otimes \V_i^+ \otimes \Sym^2 \V_i^-,\, 
\V_i^- \otimes \V_i^- \otimes \Sym^2 \V_i^+, \,
\Sym^2 \V_i^+  \otimes \Sym^2 \V_i^-. \]
Since each of them contains the trivial $\gs^i$-module, we conclude that \eqref{eq411} contains the trivial $\gs$-module.

\np
Finally, we consider the case when $\delta_{i_j} = \delta_{i_{j+1}}$ for some $1 \leq j \leq l$. (The minimality of the cycle \eqref{eq31} implies that
there are at most two such indices but we will not use this observation.) We split the roots $\nu_1, \nu_2, \ldots, \nu_l$ into two groups 
$\lambda_1, \lambda_2, \ldots, \lambda_s$ and $\mu_1, \mu_2, \ldots, \mu_t$ in the following way: If $\nu_j = s_j \delta_{i_j} - s_{j+1} \delta_{i_{j+1}}$,
then we put $\nu_j$ in the first or second group depending on whether it appears once or twice in $\nu_1, \nu_2, \ldots, \nu_l$, if $\nu_j = s_j \delta_{i_j}$, 
we put $\nu_j$ in the second group, and if $\nu_j =  2 s_j \delta_{i_j}$, we put $\nu_j$ in the first group. Set $l':= s + 2t$; note that $l' \neq l$. 

\np
From this point on the argument repeats the argument above with the following modifications:
\begin{enumerate}
\item[(i)] We consider $\Sym^{l'} (\gg^{\lambda_1} \oplus \cdots \oplus \gg^{\lambda_s} \oplus \gg^{\mu_1} \oplus \cdots \oplus \gg^{\mu_t})$ in place of 
$\Sym^l (\gg^{\lambda_1} \oplus \cdots \oplus \gg^{\lambda_s} \oplus \gg^{\mu_1} \oplus \cdots \oplus \gg^{\mu_t})$.
\item[(ii)] In the case when $\gg = D_n$ and $(\cP, \sigma)$ is of Type I, we replace $\Sym^2 \V_{i_j}^{s_j} \otimes \Sym^2 \V_{i_{j+1}}^{-s_{j+1}}$
by $ \Lambda^2 \V_{i_j}^{s_j} \otimes \Lambda^2 \V_{i_{j+1}}^{-s_{j+1}}$ in \eqref{eq412}. Correspondingly, $\V_i$  appears in one of the following groups
\[\V_i^+ \otimes \V_i^+ \otimes \V_i^- \otimes \V_i^-, \,
\V_i^+ \otimes \V_i^+ \otimes \Lambda^2 \V_i^-,\, 
\V_i^- \otimes \V_i^- \otimes \Lambda^2 \V_i^+, \,
\Lambda^2 \V_i^+  \otimes \Lambda^2 \V_i^-. \]
\end{enumerate}
Exactly as above, for $i>0$, each of the groups above contains the trivial module $\gs^i$-module. Finally, 
if $\gg = \B_n$ or $\gg = \D_n$ and $(\cP, \sigma)$ is of Type II, $\V_0$ appears in groups $\Sym^2 \V_0$ 
(one for each $\nu_j = s_j \delta_{i_j}$).  Since in these cases $\gs^0 = \B_{|\I_0|}$ or 
$\gs^0 = \D_{|\I_0|}$, $\Sym^2 \V_0$ contains the trivial $\gs^0$-module. This completes the proof.
\qed

\np
We now turn to the case that $\Ssh$ is not saturated.

\np 
\point 
{\bf Existence of $\gp_{\Msh}$ in types $\A$ and $\D$.} 
If $\gg$ is of type $\A$ there is nothing to prove since every subset $\Rsh$ is saturated and the 
statement is equivalent to the first part of this section.
The situation is the same when $\gg = \D_n$ and $(\cP, \sigma)$ is of type I. 

\np
Let $\gg = \D_n$ and let $(\cP, \sigma)$ be of type II. We will extend the proof of part ({\em a}) to this case.

\np
First we note that $-2\delta_i \in \Ssh$  and $\delta_i \in \Ssh$ imply that $(\Sym^\cdot(\Msh))^{\gs} \neq \CC$.
Indeed, $\Lambda^2 \V_i^- \oplus \V_i^+ \otimes \V_0$ is a submodule of $\Msh$ and hence we have the following inclusions of modules:

\begin{equation}\label{eq55}
\begin{array}{rcl}
\Sym^6 (\Lambda^2 \V_i^- \oplus \V_i^+ \otimes \V_0)& \subset & \Sym^\cdot(\Msh)\\
\Sym^2(\Lambda^2 \V_i^-)  \otimes \Sym^4 (\V_i^+ \otimes \V_0) & \subset & \Sym^6 (\Lambda^2 \V_i^- \oplus \V_i^+ \otimes \V_0)\\
\SS^{(2,2)} \V_i^+ \otimes \SS^{(2,2)} \V_0  \subset  \Sym^4 (\V_i^+ \otimes \V_0) & , &
\SS^{(2,2)} \V_i^-  \subset  \Sym^2(\Lambda^2 \V_i^-),
\end{array} 
\end{equation}
where $\SS^{(2,2)} \WW$ denotes the result of applying the Schur functor $\SS^{(2,2)}$ to $\WW$. 
The above inclusions along the fact that
$\SS^{(2,2)} \V_0$ contains the trivial $\gs^0$-module imply that $(\Sym^6(\Msh))^\gs \neq 0$. A symmetric argument shows that 
$2\delta_i \in \Ssh$  and $-\delta_i \in \Ssh$ imply that $(\Sym^\cdot(\Msh))^{\gs} \neq \CC$.

\np
From this point on the proof follows the proof of part ({\em a}) with the following modifications:
\RomanList
\begin{enumerate}
\item In the definition of $\prec$ we use $s_i \delta_i \prec - s_i \delta_i$ if $s_i \delta_i \in \Ssh$ {\bf or}  $2s_i \delta_i \in \Ssh$.
\item If $s_i \delta_i \prec - s_i \delta_i$, $\nu_i$ denotes the corresponding element of $\Ssh$ above; if there are two such elements, we 
set $\nu_i := s_i \delta_i$.
\item In splitting $\nu_1, \nu_2, \ldots, \nu_l$ into two groups $\lambda_1, \lambda_2, \ldots, \lambda_s$ and $\mu_1, \mu_2, \ldots, \mu_t$,
we put a root $\nu_i$ from (ii) into the first group if $\nu_i = 2s_i \delta_i$ and in the second group otherwise.
\item We consider $\Sym^{2l'} (\gg^{\lambda_1} \oplus \cdots \oplus \gg^{\lambda_s} \oplus \gg^{\mu_1} \oplus \cdots \oplus \gg^{\mu_t})$ in place of 
$\Sym^l (\gg^{\lambda_1} \oplus \cdots \oplus \gg^{\lambda_s} \oplus \gg^{\mu_1} \oplus \cdots \oplus \gg^{\mu_t})$.
\item We replace the module in \eqref{eq411} by 
$\Sym^2\gg^{\lambda_1} \otimes \cdots \otimes \Sym^2\gg^{\lambda_s} \otimes \Sym^4\gg^{\mu_1} \otimes \cdots \otimes \Sym^4 \gg^{\mu_t}$.
\end{enumerate}
\NumberList 
Using the inclusions \eqref{eq55} we conclude that $(\Sym^\cdot(\Msh))^{\gs} \neq \CC$. 
This completes the proof when $\gg = \D_n$. \qed
%This proves part ({\em b}) for $\gg$ of type $\A$ or $\D$.

\np
\point 
{\bf Examples in types $\B$ and $\C$ when $\Msh$ is not saturated.}
We will now construct examples in types $\B$ and $\C$ 
of $\gs$ and $\Ssh$ such that $(\Sym^\cdot(\Msh))^{\gs} = \CC$ and for which there does not exist a
parabolic subalgebra $\gp_{\Msh}$ of $\gg$ with reductive part $\gs$ and $\Msh \subset \gp_{\Msh}$.

\np
If $\gg = \B_n$, consider $\gs = \gs_{(\cP, \sigma)}$, where $\cP$ is the partition of type I
\[\{1,2\} \prec \{3\} \prec \{4\} \prec \cdots \prec \{n\}\]
and $\sigma(i) = 1$ is constant. Then $\gs^1 = \ggl_2$. Moreover, $\U:= \gg^{-\delta_1}$ is the $\ggl_2$-module which is the natural representation of $\gsl_2$ and
on which the identity matrix of $\ggl_2$ acts as multiplication by $-1$ and $\WW:=\gg^{2\delta_1}$ is the 
one dimensional $\ggl_2$-module on which the identity
matrix acts as multiplication by $2$. Let $\Ssh := \{-\delta_1, 2 \delta_1\}$. Then $\Msh = \U \oplus \WW$ and 
\[\Sym^k \Msh = \oplus_j \Sym^j \U \otimes \Sym^{k-j} \WW.\]
Note that $\Sym^j \U \otimes \Sym^{k-j} \WW$ is the irreducible $\gsl_2$-module of dimension $j+1$ on which the identity matrix of $\ggl_2$ acts as 
multiplication by $2k - 3j$. This proves that $(\Sym^\cdot(\Msh))^\gs = \CC$ but there is no 
parabolic subalgebra $\gp_{\Msh}$ of $\gg$ with reductive part $\gs$ such that $\Msh \subset \gp_{\Msh}$.

\np
If $\gg = \C_n$, consider $\gs = \gs_{(\cP, \sigma)}$, where $\cP$ is the partition of type II
\[\{1\} \prec \{2\} \prec \{3\} \prec \cdots \prec \{n\}\]
and $\sigma(i) = 1$ is constant. Then $\gs^0 = \C_1 \cong \gsl_2$ and $\gs^1 = \ggl_1$, i.e. $\gs^0 \oplus \gs^1 \cong \ggl_2$. Moreover, setting
$\U:= \gg^{-\delta_1}$ and  $\WW:=\gg^{2\delta_1}$, we arrive at exactly the same situation as in the case 
$\gg = \B_n$ above.  \qed

%\section{Proof of Theorem \ref{thm:classical} in the case when $\gg$ is exceptional.}
\section{Proof of the Main Theorem when $\gg$ is exceptional.}
\label{sec:proof-exceptional}

\np 
\point
First we recall some standard notation following the conventions in \cite{Bo}. If $\gg$ is a simple Lie algebra of rank $n$
we label the simple roots of $\gg$ as $\alpha_1, \ldots, \alpha_n$ as in \cite{Bo}. The fundamental dominant weight of $\gg$ are denoted by
$\omega_1, \ldots, \omega_n$. If $-\alpha_0$ is the highest root, then
$\alpha_0, \alpha_1, \ldots, \alpha_n$ label the extended Dynkin diagram of $\gg$.

\np 
\point 
{\bf Existence of $\gp_{\Msh}$ in type $\G_2$ when $\Ssh$ is saturated.} 
Let $\gg = \G_2$. Let
$\Ssh$ be a saturated subset of $\Rsh$
and let $\Msh = \oplus_{\nu \in \Ssh} \gg^\nu$. If 
$(\Sym^\cdot(\Msh))^{\gs} = \CC$, then there exists a
parabolic subalgebra $\gp_{\Msh}$ of $\gg$ with reductive part $\gs$ such that $\Msh \subset \gp_{\Msh}$.
Indeed, if $\gs$ is a proper subalgebra of $\gg$ which is not equal to
$\gh$, then all elements of $\Rsh$ are proportional and there is nothing to prove. If $\gs = \gh$, then the spaces $\gg^\nu$ are just
the root spaces of $\gg$ which are one dimensional and again the statement is clear. 
\qed

\np
\point 
{\bf Example in type $\G_2$ when $\Ssh$ is not saturated.} 
On the other hand, let $\gs \cong \ggl_2 \subset \gg$ be the parabolic subalgebra of $\gg$ with roots $\pm \alpha_2$. 
Then $\Rsh = \{\pm \delta, \pm 2 \delta, \pm 3 \delta\}$. Moreover, $\gg^{\pm k \delta}$ is the irreducible $\gs$-module 
of dimension $2, 1$, or $2$ (corresponding to $k = 1, 2$, or $3$) on which a fixed element in the centre of $\gs$
acts as multiplication by $\pm k$. Then, for $\Ssh = \{-\delta, 2 \delta\}$, setting $\U:= \gg^{-\delta}$ and  $\WW:=\gg^{2\delta}$, 
we arrive at exactly the same situation as at the end of Section 2 above. In particular,
$(\Sym^\cdot(\Msh))^\gs = \CC$ but there is
no parabolic subalgebra $\gp_{\Msh}$ of $\gg$ with reductive part $\gs$ such that $\Msh \subset \gp_{\Msh}$.
\qed

\np
\point 
{\bf Examples in types $\F_{\!4}$, $\E_6$, $\E_7$, and $\E_8$ with $\Ssh$ saturated.} 
Let $\gg = \F_{\!4}, \E_6, \E_7$, or $\E_8$.  We will construct a saturated set $\Ssh$ such that 
$(\Sym^\cdot(\Msh))^{\gs} = \CC$ but there is no
parabolic subalgebra $\gp_{\Msh}$ of $\gg$ with reductive part $\gs$ such that $\Msh \subset \gp_{\Msh}$.

\np
Denote the rank of $\gg$ by $n$. Consider the extended Dynkin diagram of $\gg$. Removing the node
connected to the root $\alpha_0$ we obtain the Dynkin diagram 
of a semisimple subalgebra $\gm \oplus \gc$ of $\gg$ of rank $n$ where $\gm \cong \A_1$ is the subalgebra of $\gg$ with
roots $\{\pm \alpha_0\}$ and $\gc$ is the subalgebra if $\gg$ with simple roots obtained from the simple roots of $\gg$ after removing 
the one adjacent to $\alpha_0$. More precisely, we remove the roots 
$\alpha_1, \alpha_2, \alpha_1, \alpha_8$ when $\gg = \F_4, \E_6, \E_7, \E_8$ respectively. The respective subalgebras
$\gc \subset \gg$ are isomorphic to 
$\gc \cong \C_3, \A_5, \D_6$, or $\E_7$ respectively. As an $\gm$--module $\gg$ decomposes as

\begin{equation} \label{eq41}
\gg = (\Ad_\gm \otimes \tr_\gc) \oplus (\tr_\gm \otimes \Ad_\gc) \oplus (\V \otimes \U),
\end{equation}
where $\Ad_\gm$ and $\Ad_\gc$ are the adjoint modules of $\gm$ and $\gc$ respectively; $\tr_\gm$ and $\tr_\gc$ ---the respective trivial modules;
$\V$ is the natural $\gm \cong \A_1$--module; and $\U$ is the $\gc$--module whose highest weight is the fundamental weight of $\gc$ corresponding
to the simple root of $\gc$ linked to the removed node of the extended Dynkin diagram of $\gg$. In fact, for $\gg = \F_4, \E_6, \E_7, \E_8$,
the highest weight of $\gc$ is $\omega_3, \omega_3, \omega_6, \omega_7$ respectively.
Here the weights of $\U$ are given according to the labeling conventions of $\gc$. For example, if $\beta_1, \beta_2, \beta_3$ are the simple roots
of $\gc = \C_3$ in the case when $\gg = \F_4$, we have $\beta_1 = \alpha_4$, $\beta_2 = \alpha_3$, and $\beta_3 = \alpha_2$. 

\np
Set $\gs = \gm + \gh$. From the construction of $\gs$ we conclude that $\gt = \gh_\gc$, the Cartan subalgebra of $\gc$. Furthermore,
 (\ref{eq41}) implies $\Rsh = \Delta_\gc \cup \supp \U$, where $\supp \U$ denotes the set of weights 
 of $\U$ and, for $\nu \in \Rsh$ the 
 $\gs = \gm \oplus \gh_\gc$--module $\gg^\nu$ is given by
 $$
 \gg^\nu \cong \begin{cases}
 \tr_\gm \otimes \nu & {\text { if }} \nu \in \Delta_\gc\\
 \V \otimes \nu & {\text { if }} \nu \in \supp \U.
 \end{cases}
 $$
Let $\omega$ be the highest weight of  $\U$ and write $\omega= q_1 \beta_1 + \cdots + q_{n-1} \beta_{n-1}$ where 
$\beta_1, \ldots, \beta_{n-1}$ are the simple roots of $\gc$ and $q_i \in \QQ_+$. Set
$\Ssh = \{ -\omega, \beta_1, \ldots, \beta_{n-1}\}$. Then $\Msh=\gg^\omega \oplus  (\opp_{i=1}^{n-1} \gg^{\beta_i})$ and
\[ \Sym^k \Msh = \bigoplus_{j+ i_1 + \cdots + i_{n-1} = k} \Sym^j \gg^{-\omega} \otimes \Sym^{i_1} \gg^{\beta_1} \otimes \cdots \otimes \Sym^{i_{n-1}} \gg^{\beta_{n-1}}.\]
Moreover, $\Sym^j \gg^{-\omega} \otimes \Sym^{i_1} \gg^{\beta_1} \otimes \cdots \otimes \Sym^{i_{n-1}} \gg^{\beta_{n-1}}$ is an irreducible $\gm$--module
which is not trivial unless $j = 0$ and on which $\gh_\gc$ acts via $-j \omega + i_1 \beta_1 + \cdots + i_{n-1} \beta_{n_1}$. This implies that, for $k >0$,
$(\Sym^k \Msh)^\gs = 0$ and hence $(\Sym^\cdot \Msh)^\gs = \CC$. 
On the other hand, the equation $\omega= q_1 \beta_1 + \cdots + q_{n-1} \beta_{n-1}$ 
implies that there is no
parabolic subalgebra $\gp_{\Msh}$ of $\gg$ with reductive part $\gs$ such that $\Msh \subset \gp_{\Msh}$.
\qed

%\bigskip
%\bigskip
%\np
%{\small
%Department of Mathematics and Statistics, \\
%Queen's University, Kingston, Ontario, Canada, K7L 3N6.\\
%{\em E-mail address}: {\tt dimitrov@mast.queensu.ca}, {\tt mikeroth@mast.queensu.ca}}

\end{document}